\documentclass[a4paper]{amsart}
\usepackage{amssymb, amsmath, latexsym}
\usepackage[all]{xy}
\begin{document}
\newcommand{\pn}{\par \noindent}
\newcommand{\pmn}{\par \medskip \noindent}
\newcommand{\pbn}{\par \bigskip \noindent}
\newtheorem{De}{Definition}[section]
\newtheorem{Th}[De]{Theorem}
\newtheorem{Pro}[De]{Proposition}
\newtheorem{Fac}[De]{Fact}
\newtheorem{Le}[De]{Lemma}
\newtheorem{Co}[De]{Corollary}
\theoremstyle{definition}
\newtheorem{Ob}[De]{Observation}
\newtheorem{Rem}[De]{Remark}
\newtheorem{Exa}[De]{Example}
\newcommand{\Def}[1]{\begin{De}#1\end{De}}
\newcommand{\Thm}[1]{\begin{Th}#1\end{Th}}
\newcommand{\Prop}[1]{\begin{Pro}#1\end{Pro}}
\newcommand{\Fact}[1]{\begin{Fac}#1\end{Fac}}
\newcommand{\Obs}[1]{\begin{Ob}#1\end{Ob}}
\newcommand{\Ex}[1]{\begin{Exa}#1\end{Exa}}
\newcommand{\Lem}[1]{\begin{Le}#1\end{Le}}
\newcommand{\Cor}[1]{\begin{Co}#1\end{Co}}
\newcommand{\Rek}[1]{\begin{Rem}#1\end{Rem}}
\newcommand{\ra}{\rightarrow}
\newcommand{\lra}{\longrightarrow}
\def\leq{\leqslant}
\def\geq{\geqslant}
\newcommand{\ga}{\Gamma}
\newcommand{\Ga}{\Gamma}
\def\HG{H\!\varGamma}
\def\HGC{H\!\varGamma\! C}
\def\HH{\mathit{H\!H}}
\def\HC{\mathit{H\!C}}
\newcommand{\sis}{{\cal S}_*}
\def\Sym{\operatorname{Sym}}
\def\Tym{\operatorname{Tym}}
\newcommand{\K}{\mathbb{K}}
\newcommand{\Q}{\mathbb{Q}}
\newcommand{\F}{\mathbb{F}}
\newcommand{\D}{\mathbb{D}}
\newcommand{\Z}{\mathbb{Z}}
\newcommand{\s}{\mathbb{S}}
\newcommand{\taq}{{\HG}}
\newcommand{\aq}{{\sf AQ}}
\newcommand{\tot}{{\sf Tot}}
\newcommand{\tor}{{\sf Tor}}
\newcommand{\nat}{\sf Nat}
\newcommand{\Untor}{{\sf Untor}}
\newcommand{\cotor}{{\sf Cotor}}
\newcommand{\stab}{{\sf stab}}
\newcommand{\einf}{E_{\infty}}
\newcommand{\hot}{\hat{\otimes}}
\newcommand{\homo}{{\rm Hom}}
\newcommand{\call}{{\mathcal L}}
\newcommand{\calc}{{\mathcal C}}
\newcommand{\calh}{{\mathcal H}}
\newcommand{\cala}{{\mathcal A}}
\newcommand{\calf}{{\mathcal F}}
\newcommand{\unn}{\underline{n}}
\newcommand{\unm}{\underline{m}}
\newcommand{\ot}{\overline{t}}
\newcommand{\ie}{\emph{i.e.}}
\newcommand{\eg}{\emph{e.g.}}
\title[The collapse of the periodicity sequence]{The collapse of the periodicity sequence
in the stable range}
\date{\today}
\author{Birgit Richter}
\address{Fachbereich Mathematik der Universit\"at Hamburg,
Bundesstrasse 55, 20146 Hamburg, Germany}
\email{richter@math.uni-hamburg.de}
\begin{abstract}
The stabilization of Hochschild homology of commutative algebras
is Gamma homology. We describe a cyclic variant of Gamma homology
and prove that the associated analogue of Connes' periodicity sequence
becomes almost
trivial, because the cyclic version coincides with the ordinary
version from homological degree two on. We offer an alternative
explanation for this by proving that the $B$-operator followed by
the stabilization map is trivial from degree one on.
\end{abstract}

\maketitle
\section{Introduction}

Given a commutative algebra $A$ there are several homology theories
available that can help to understand $A$. Hochschild homology and
cyclic homology of $A$ are related by Connes' periodicity sequence
$$
\ldots \ra \HH_n(A) \stackrel{I}{\lra} \HC_{n}(A)  \stackrel{S}{\lra}
\HC_{n-2}(A)  \stackrel{B}{\lra} \HH_{n-1}(A) \ra
\ldots$$
which is a good means for comparing Hochschild homology with its
cyclic variant.

Using the commutativity of $A$ we could consider Andr\'e-Quillen
homology as well. Viewing $A$ as an $E_\infty$-algebra with
trivial homotopies for commutativity allows us to consider
Andr\'e-Quillen homology  in the category of differential graded
$E_\infty$-algebras as defined by Mike Mandell \cite{M}. This
homology theory coincides with Alan Robinson's Gamma homology
\cite{Ro, BR} which in turn can be interpreted as
stabilization of Hochschild homology of $A$ by \cite[Theorem
1]{PR}.

This homology theory has the feature that it coincides with
Andr\'e-Quillen homology for $\mathbb{Q}$-algebras \cite[Theorem 6.4]{RoWh}.

The Hodge decomposition for Hochschild homology for commutative algebras whose
base ring contain the rationals splits Andr\'e-Quillen homology off
as the first summand $\HH_*^{(1)} \cong AQ_{*-1}$ in the
decomposition (see \cite{GS,L2,NS}). Cyclic homology splits
similarly and from degree three on the first summand $\HC_*^{(1)}$
of that decomposition is again Andr\'e-Quillen homology, $AQ_{*-1}$.
It is known  that the periodicity sequence passes to a sequence for
the decomposition summands \cite{L2,NS},

$$ \ldots \ra \HH_n^{(i)}(A) \stackrel{I}{\lra} \HC_n^{(i)}(A) \stackrel{S}{\lra}
\HC^{(i-1)}_{n-2}(A) \stackrel{B}{\lra} \HH_{n-1}^{(i)}(A) \ra
\ldots$$

Therefore rationally the periodicity sequence collapses in higher
degrees for the first decomposition summand, because there the map $I$
becomes an isomorphism. One could guess that this is a defect of
working over the rationals, but we will show in the course of this
paper that this is not the case.

Robinson and Whitehouse proposed a cyclic variant of Gamma
homology of differential graded $E_\infty$-algebras over a cyclic
$E_\infty$-operad in \cite{RoWh}.

Motivated by the application of Gamma homology to obstruction theories for
$E_\infty$ ring structures on ring spectra (see \cite{Ro}) we construct
a cyclic variant of Gamma homology in the restricted case of commutative
algebras which arises naturally from the interpretation of Gamma homology
as stable homotopy of certain $\Gamma$-modules.

The aim of this paper is to deduce a periodicity sequence for Gamma
homology and its cyclic version. It turns out, however,
(\emph{cf.} \ref{cor:iso}) that
cyclic Gamma homology coincides with usual Gamma homology from
homological degree two on; hence this
sequence collapses. We can explicitly describe (see
Propositions \ref{prop:deg0} and \ref{prop:lowdeg}) cyclic Gamma
homology in  small degrees in terms of ordinary cyclic homology and
deRham cohomology.

As one explanation for this behavior we show (Theorem
\ref{thm:collapse})  that the composition of the stabilization map
with the $B$-operator is trivial in degrees bigger than zero. For
large enough degrees the stable and unstable periodicity sequence
are related as follows:
$$ \xymatrix{
{\ldots} \ar[r] & {\scriptstyle \HC_{n+2}(A)} \ar[r]^{S} &
{\scriptstyle \HC_n(A)} \ar[rd]^{0} \ar[r]^{B} &
{\scriptstyle \HH_{n+1}(A)} \ar[r]^{I}  \ar[d]^{\stab}& {\scriptstyle
  \HC_{n+1}(A)}  \ar[r]^{S} & {\ldots} \\
{\ldots} \ar[r] & {\scriptstyle \HGC_{n+1}(A)} \ar[r] & {\scriptstyle 0
} \ar[r] &
{\scriptstyle \HG_n(A)} \ar[r]^{\cong} & {\scriptstyle \HGC_n(A)}
\ar[r] &  {\ldots}
}$$

Our results must disappoint everybody who hoped that a cyclic version
of Gamma homology would help to calculate Gamma homology and probably
identify obstruction classes in the setting of \cite{Ro}. However, we
clarify the r{\^o}le Gamma homology plays as a stabilization of Hochschild
homology.

The proofs use the extension of the definitions of cyclic,
Hochschild and Gamma homology to functor categories. We recall the
necessary prerequisites from \cite{L,P,PR}.


Lars Hesselholt proved an analogous phenomenon in the setting of
topological Hochschild homology. Fix an arbitrary prime $p$.
In \cite{H} he showed
that the equivalence between stable $K$-theory and topological
Hochschild homology is reflected in an equivalence between the
$p$-completions of the stabilization of topological cyclic homology
and $p$-completed topological Hochschild homology.

\section{The category $\calf$ and $\calf$-modules } \label{calf}
We recall the definition of cyclic homology from \cite[\S 6]{L}
(see also \cite[\S 3]{P}). Let $\calf$ denote the skeleton of the
category of  finite unpointed sets and let $\unn$ be the object
$\{0,\ldots,n\}$ in  $\calf$. We call functors from $\calf$ to the
category of  $k$-modules $\calf$-modules. Here $k$ is an arbitrary
commutative ring with  unit. For a set $S$ we denote by  $k[S]$ the free
$k$-module generated by $S$.  

The projective generators for the category of
$\calf$-modules are the  functors $\calf^n$ given by
$$ \calf^n(\unm) := k[\calf(\unn, \unm)]; $$
whereas the category of contravariant functors from $\calf$ to $k$-modules has
the family $\calf_n$ with
$$ \calf_n(\unm) := k[\calf(\unm, \unn)] $$
as generators.

For two $\calf$-modules $F$ and $F'$ let $F \otimes F'$ be the
pointwise tensor product of $F$ and $F'$, i.e., $F \otimes
F'(\unn) = F(\unn) \otimes F'(\unn)$. As a map in $\calf$ from the
object $\underline{0}$ to an object $\unm$ just picks an arbitrary
element, one obtains that $(\calf^0)^{\otimes n}
 \cong \calf^{n-1}$. The functor $\calf^0$ is an analog of the functor $L$ from
 \cite{PR} in the unpointed setting and the tensor powers
$(\calf^0)^{\otimes n} \cong \calf^{n-1}$ for
$n > 1$ correspond to $L^{\otimes n}$.

Given a unital commutative 
$k$-algebra $A$, the $\calf$-module which gives rise to cyclic
homology of $A$ is the functor $\call(A)$ that sends $\unn$ to
$A^{\otimes
  n+1}$.  A map $f\colon \unn
\ra \unm$ induces $f_*\colon \call(A)(\unn) \ra \call(A)(\unm)$ via
$$f_*(a_0 \otimes \ldots \otimes a_n) = b_0 \otimes \ldots \otimes
b_m, \text{ with } b_i = \prod_{j \in f^{-1}(i)} a_j.$$
Here we set $b_i = 1$ if the preimage of $i$ is empty.

For any $\calf$-module $F$, cyclic homology of $F$, $\HC_*(F)$, can be defined
\cite[3.4]{P} as the homology of the total complex associated to the bicomplex

$$ \xymatrix{
{\vdots}\ar[d]^{b} & {\vdots} \ar[d]^{b} & {\vdots} \ar[d]^{b} \\
{F(\underline{2})}
\ar[d]^{b} & {F(\underline{1})} \ar[l]^{B} \ar[d]^{b} &
{F(\underline{0})}  \ar[l]^{B}\\
{F(\underline{1})} \ar[d]^{b} & {F(\underline{0})} \ar[l]^{B}\\
{F(\underline{0})} } $$
In particular, cyclic homology of $A$, $\HC_*(A)$, is the homology of
this total complex applied to the functor $\call(A)$. We recall the
definition of $b$ in \eqref{eqn:b} and the one of $B$ in Definition
\ref{defn:B}.

\section{The relationship to the category $\ga$}
Let $\ga$ be the skeleton of the category of pointed finite sets
and let $[n]$ be the object $[n] = \{0,\ldots,n\}$ with $0$ as
basepoint. The projective generators of the category of
$\ga$-modules are the functors $\ga^n$ given by
$$ \ga^n[m] = k[\ga([n],[m])]. $$
There is a natural forgetful functor $\mu\colon \ga \ra \calf$ and the left
adjoint to $\mu$,
$\nu\colon \calf \ra \ga$,  which adds an extra basepoint $\nu(\unm) =
[m+1]$. Pulling  back with these functors transforms $\ga$-modules
into  $\calf$-modules and vice  versa:
$$ \xymatrix{ {\mu^*\colon \calf-{\rm modules}} \ar@<-0.5ex>[r] &
{\ga-{\rm modules} \colon\nu^*}   \ar@<-0.5ex>[l]
}$$
In \cite[Proposition 3.3]{P} Pirashvili shows that
$$ \tor_*^\calf(\nu^*F,G) \cong \tor_*^\ga(F,\mu^*G).$$
\Lem{The functor $\calf^n$ pulled back along $\mu$ is
isomorphic to  $\ga^{n+1}$. \label{lemma}
}
\begin{proof}
We first show that the $\ga$-module $\mu^*(\calf^0)$ is
isomorphic to $\ga^1$: On every object $[n]$ we obtain that
$$\mu^*(\calf^0)[n]= \calf^0(\unn) = k[\calf(\underline{0}, \unn)] \cong
k^{n+1} $$
because the value of a function $f \in \calf(\underline{0}, \unn)$ on
$0$ can be an arbitrary element $i \in \unn$. The $\ga$-module $\ga^1$
has the  same value
on $[n]$, because a function $g \in \ga([1],[n])$ has an arbitrary value on $1$
 but sends zero to zero. As we just allow pointed maps,  the two functors are
isomorphic.

The general case easily follows  by direct considerations or by using the
decompositions of $\calf^n $ and $\ga^{n+1}$  as $(n+1)$-fold
tensor products  $\calf^n
\cong (\calf^0)^{\otimes n+1}$ and $\ga^{n+1} \cong (\ga^1)^{\otimes
  n+1}$.
\end{proof}
Recall, that Hochschild homology of a $\ga$-module $G$, $\HH_*(G)$,
can be defined as the homology of the complex
\begin{equation} \label{eqn:b}
 G[0] \xleftarrow{b} G[1] \xleftarrow{b} \cdots
\end{equation}
where $b = \sum_{i=0}^n G(d_i)$ and $d_i$ is the map of pointed sets
that for $i<n$ sends $i$ and $i+1$ to $i$ and is bijective and order
preserving on the other values in $[n]$. The last map, $d_n$, maps $0$
and  $n$ to $0$ and is the identity for all other elements of $[n]$.

Later, we will need the following auxiliary result.

\Lem{ \label{lem:deg}
Hochschild homology of a Gamma module $G$ is isomorphic to the
homology of the normalized complex which consist of
$G[n]/D_n$ in chain degree $n$ where $D_n \subset G[n]$ consists of
all elements
of the form $(s_i)_*F[n-1]$ where $s_i$ is the order preserving injection from
$[n-1]$ to $[n]$ which misses $i$.
}
\begin{proof}
This result just uses the standard fact that the Hochschild complex
is the chain complex associated to a simplicial $k$-module and the
elements in $D_n$ correspond to the degenerate elements; therefore
the complex $D_*$ is acyclic.
\end{proof}
A similar result applies to cyclic homology of $\calf$-modules.

Let $\mathbb{S}^1 = \Delta^1/\partial \Delta^1$ denote the
standard model of the simplicial $1$-sphere.  Recall from
\cite{L,P} that Hochschild homology of a commutative unital
$k$-algebra $A$, $\HH_*(A)$,  coincides with the homotopy groups of the
simplicial $k$-module $\mu^* \call(A)(\mathbb{S}^1)$. Here, we
evaluate $\mu^* \call(A)$ degreewise. More general, Hochschild
homology of any $\ga$-module $G$ coincides with $\pi_*G(\mathbb{S}^1)$.

\section{Gamma homology and its cyclic version}
Let $t$ be the contravariant functor from $\ga$ to $k$-modules
which is defined as
$$ t[n] = \homo_{\sf Sets_*}([n],k) $$
where ${\sf Sets_*}$ denotes the category of pointed sets.
Pirashvili and the author proved in \cite{PR} that Gamma homology
of any $\ga$-module $G$, $\HG_*(G)$, is isomorphic to
$\tor_*^\ga(t, G)$. In particular, Gamma homology of the algebra
$A$, $\HG_*(A)$ is isomorphic to $\tor_*^\ga(t, \mu^*\call(A))$.

For a cyclic variant of Gamma homology, we have to transform the
functor $t$ into a contravariant $\calf$-module. Choosing $\nu^*t$
does this, but it inserts an extra basepoint. Killing the value on
an additional point amounts to define the $\calf$-module $\ot$ by
the following exact sequence:
$$ 0 \lra \calf_0 \lra \nu^*(t) \lra \ot \lra 0. $$
The transformation from $\calf_0$ to $\nu^*t$ is given by sending a
scalar multiple $\lambda f$ of a map $f\colon \unn \ra \underline{0}$ to
the function in $\homo_{\sf Sets_*}([n+1],k)$ which sends the points
$1,\ldots,n+1$ to $\lambda$.

\Prop{
On the family of projective generators $\left(\calf^n\right)_{n \ge 0}$
the torsion groups with respect to $\ot$ are as follows: \label{proj}
$$\tor_*^\calf(\ot,\calf^n) \cong \left\{ \begin{array}{rl}
0 & {\rm for }  \quad * > 0\\
k^n &  {\rm for }  \quad * = 0
\end{array} \right.$$
}
\begin{proof}
It is clear that the torsion groups vanish in positive degrees
because the functors $\calf^n$ are projective. We have to
 prove the claim in  degree zero, but the tensor products
in question are easy to calculate:
$$ \ot \otimes_\calf \calf^n \cong \ot(\unn) \cong k^n.$$
\end{proof}
\Def{We call the group $\tor_n^\calf(\overline{t}, F)$ the \emph{$n$th
cyclic Gamma homology group of the $\calf$-module $F$} and
denote it by $\HGC_n(F)$.}

\Rek{We will see in \ref{prop:deg0} that cyclic Gamma homology of
an algebra $A$ in degree zero behaves analogously to usual Gamma
homology whose value in homological degree zero gives Hochschild
homology of degree one.}

As the functor $\calf_0$ is projective, the calculation in
proposition \ref{proj} allows us to  draw the following conclusion.
\Cor{\label{cor:iso} Cyclic Gamma homology of any $\calf$-module $F$
coincides with Gamma homology of the induced Gamma module $\mu^*(F)$ in
degrees higher than $1$, \ie,
$$ \HGC_*(F) = \tor_*^\calf(\overline{t}, F) \cong
\tor_*^\ga(t, \mu^*F) \cong \HG_*(\mu^*(F)) \quad \forall * >
1.$$}

In low degrees the difference between cyclic and ordinary Gamma
homology is measured by the following exact sequence:
$$ 0 \ra \tor_1^\calf(\nu^*t, F) \ra \tor_1^\calf(\overline{t}, F)
\stackrel{\delta}{\lra}
F(\underline{0}) \ra \nu^*t \otimes_{\calf} F \ra \overline{t}
\otimes_{\calf} F \ra 0 $$
which is nothing but
$$ 0 \ra \HG_1(\mu^*F) \ra \HGC_1(F) \ra \mu^*F(\underline{0})
\stackrel{\delta}{\lra}
\HG_0(\mu^*F) \ra \HGC_0(F) \ra 0. $$
We will obtain more explicit descriptions in  the algebraic case in
the next
section.
\section{The $B$ operator}
In the unstable situation
there is a map $B$ which connects cyclic homology and Hochschild homology and
which gives rise to Connes' important periodicity sequence
$$ \cdots \lra  \HH_n(A) \stackrel{I}{\lra} \HC_n(A) \stackrel{S}{\lra}
\HC_{n-2}(A) \stackrel{B}{\lra}  \HH_{n-1}(A) \lra \cdots $$
In low degrees the map $B$ sends the zeroth cyclic homology of a $k$-algebra
$A$ which is nothing but $A$ again to the first Hochschild homology group of
$A$ which consists of the module of K\"ahler differentials $\Omega^1_{A|k}$ and
the map is given by $B(a) = da$. If we consider the first nontrivial parts in
the long exact sequence of Tor-groups as above, arising from the short exact
sequence $0 \ra \calf_0 \ra \nu^*t \ra \overline{t} \ra 0$ then, for the
functor $\call(A)$, we obtain
$$ \cdots \ra A \ra  \nu^*t \otimes_{\calf} \call(A) \ra \overline{t}
\otimes_{\calf} \call(A) \ra 0 $$ and $\nu^*t \otimes_{\calf}
\call(A)$ is isomorphic to the zeroth Gamma homology group of $A$
which is the module of K\"ahler differentials. The map is induced
by the natural transformation from $\calf_0$ to $\nu^*t$. The aim of
this section is to prove that this map is given by the $B$-map.

Let us recall the general definition of the $B$-map for cyclic and
Hochschild homology of functors. The $B$-map from cyclic homology to
Hochschild homology can be
viewed as a map from the $n$th generator $\calf_n$ to the
$(n+1)$st in the following manner:

\Def{\label{defn:B}
Let $\tau$ be the
generator of the cyclic group on $n+1$ (resp $n+2$) elements and
let $s$ be the map of finite sets which sends $i$ to $i+1$. Then
the $B$-map is defined as a map $B\colon \calf_{n} \ra \calf_{n+1}$. On
a generator $f\colon \unm \ra \unn$ it is $B(f) := (-1)^n(1-\tau) \circ s
\circ N \circ f$ where $N$ is the norm map $N = \sum_{i=1}^{n+1}
(-1)^i \tau^i$. }

On the part $F(\unn) \cong \calf_n \otimes_\calf F$ of the complex
for cyclic homology of $F$ this induces the usual $B$-map known
from the algebraic case $F = \call(A)$, for a commutative algebra
$A$. By the very definition of the map it is clear that it is
well-defined on the tensor product.

In our situation we apply the $B$-map to
the first column of the double complex for cyclic homology of $F$
$$ \xymatrix{
{\vdots}\ar[d]^{b} & {\vdots} \ar[d]^{b} & {\vdots} \ar[d]^{b} \\
{F(\underline{2})}
\ar[d]^{b} & {F(\underline{1})} \ar[l]^{B} \ar[d]^{b} &
{F(\underline{0})} \ar[l]^{B}\\
{F(\underline{1})} \ar[d]^{b} & {F(\underline{0})} \ar[l]^{B}\\
{F(\underline{0})} } $$ and send all other columns to zero.  In
\cite[3.2]{P} it is shown that $\nu^*\ga_n \cong \calf_n$. Using
this we obtain an  isomorphism $F(\unn) \cong \calf_n \otimes_\calf
F \cong \nu^*\Gamma_n \otimes_\calf F \cong \Gamma_n \otimes_\Gamma
\mu^*F$ and see that $B$ gives rise to a map from the  total complex
for cyclic homology of $F$ to the complex for Hochschild homology of
$\mu^*(F)$.

A verbatim translation of the  proof for (\cite[2.5.10, 2.1]{L})
in the case of a cyclic module to our setting gives the following result:

\Lem{The map $B$ is a map of chain complexes and therefore induces a
  map from $\HC_*(F)$ to $\HH_{*+1}(\mu^*F)$.}
\Rek{In degree zero, the $B$-map from $\calf_0$ to $\calf_1$ applied
  to an $f \in \calf_0(\unn)$ reduces
  to $(1-\tau) \circ s \circ f$}

We should first make sure that cyclic Gamma homology has the right value in
homological dimension zero.
\Prop{\label{prop:deg0}
Cyclic Gamma homology in degree zero is isomorphic to cyclic homology
in degree one. In particular, $\HGC_0(A) \cong \HC_1(A)$.  }
\begin{proof}
The cokernel of the map $F(\underline{0}) \ra \nu^*(t) \otimes_\calf
F$ can be determined by a map from $F(\underline{0})$ to
$F(\underline{1})$: 
similar to the beginning of the resolution $\ldots \ra \ga_2 \ra
\ga_1$ of $t$, the exactness of $\nu^*$ turns this into a resolution
$\ldots \ra \calf_1 \ra \calf_0$ of $\nu^*t$. The projectivity of
$\calf_0$ therefore gives us a lift $F(\underline{0}) \ra
F(\underline{1})$. In this lift only one summand of the $B$-map
arises: instead of the sum $(1-\tau) \circ s \circ f$ a generator $f
\in \calf_0(\unn)$ is sent to $s \circ f$. But Hochschild and cyclic
homology coincide with their normalized version (see Lemma
\ref{lem:deg}) and the second summand $\tau \circ s \circ f$ has an
image in the degenerate part.

In the algebraic case,
this lift induces a map from $A$ to $A
\otimes A$ which sends $a$ to $1 \otimes a$. The projection to the
K\"ahler differentials is then just the map $a \mapsto d(a)$ which
is the same as the $B$-map in this dimension.
\end{proof}

Cyclic Gamma homology in dimension one can be explicitly described as
well.  In small degrees our $\mathrm{Tor}$-exact sequence looks as follows:
$$ 0 \ra \HG_1(\mu^*F) \ra \HGC_1(F) \stackrel{\delta}{\lra} F(\underline{0})
\stackrel{B}{\lra} \HH_1(F). $$
Therefore we obtain the following.
\Prop{\label{prop:lowdeg}
The difference between cyclic Gamma homology and ordinary Gamma homology
in degree one is measured by the kernel of the $B$-map.}
In the
case of the functor $\call(A)$ the  exact sequence is
$$ 0 \ra \HG_1(A) \ra \HGC_1(A) \stackrel{\delta}{\lra} A
\stackrel{d}{\lra} \Omega^1_{A|k}. $$
Thus in degree one the difference between Gamma homology and its
cyclic version  is measured be the zeroth deRham cohomology of
$A$. For instance, if $A$ is \'etale, then $\HG_1(A)=0=
\Omega^1_{A|k}$ and therefore $\HGC_1(A) \cong A$.

The above calculations in small dimensions suggest that one should
view the sequence of Tor-groups coming from the sequence $0 \ra
\calf_0 \ra \nu^*t \ra \overline{t} \ra 0$ as the stable version
of the periodicity sequence. In the algebraic case the two sequences
are nicely related in the following way.
$$\xymatrix {
{\scriptscriptstyle \HC_1(A)} \ar[r]^{\scriptstyle B} &
{\scriptscriptstyle \HH_2(A)}\ar[r]^{\scriptstyle I} \ar[d]^{\stab}&
{\scriptscriptstyle  \HC_2(A)} \ar[r]^{\scriptstyle S} &
{\scriptscriptstyle  \HC_0(A)} \ar[r]^{\scriptstyle B} \ar@{=}[d]&
{\scriptscriptstyle  \HH_1(A)} \ar[r] \ar[d]^{\stab}_{\cong}&
{\scriptscriptstyle  \HC_1(A)} \ar[r]  \ar[d]_{\cong}& {\scriptscriptstyle 0}\\
{\scriptscriptstyle 0} \ar[r] &{\scriptscriptstyle \HG_1(A)} \ar[r]
&{\scriptscriptstyle  \HGC_1(A)} \ar[r]^{\delta} & {\scriptscriptstyle A}
\ar[r]^{\scriptstyle B} & {\scriptscriptstyle \HG_0(A)} \ar[r] &
{\scriptscriptstyle  \HGC_0(A)} \ar[r] &{\scriptscriptstyle 0}
}
$$

But in higher dimensions the transformation $I$ from the periodicity
sequence becomes an isomorphism. The term $\calf_0 \otimes_\calf
F\cong F(\underline{0})$ plays the role of cyclic Gamma homology in
dimension $-1$.

\section{Triviality of the $B$-map after stabilization} We will prove a
result which we like to think of as an explanation of the collapsing
of the periodicity sequence in the stable world: of course one could
say that the isomorphism of cyclic and ordinary Gamma homology in
dimensions different from zero and one explains this phenomenon, but
we would like to relate the unstable periodicity sequence to the
stable one by an explicit stabilization process.

Similar to the $B$-map, we define a stabilization map $\stab\colon 
\HH_{n+1}(G) \ra \HG_n(G)$ for a $\ga$-module $G$ on the
corresponding generator $\ga_{n+1}.$ A Gamma module $G$
is called \emph{reduced} if $G[0] =0$. As we have a unique pair of maps
$[0] \ra [n] \ra [0]$ for every $[n]$ we can split any Gamma module $G$ as
$G \cong G[0] \oplus G'$ such that  $G'$ is reduced.

Gamma homology of  a reduced functor $G'$ has a description as the
homology of the cubical construction
$Q_*(G')$ of the functor $G$ (see \cite[Theorem 4.5]{Ri} and
\cite[Theorem 1]{PR}) and $Q_*(G')$ is a tensor
product $SQ_* \otimes_\ga G'$ where $SQ_*$ is an analog of the
cubical construction of Eilenberg and MacLane on pointed sets. Gamma
homology of an unreduced functor $G$, \eg, $\call(A)$, is then just given
as
$$ \HG_*(G) = \begin{cases}
G[0] & \text{ if } *=0 \\
H_*(Q_*(G')) & \text{ if } * > 0.
\end{cases}
$$
This particular cubical model for a chain complex for Gamma homology will
give us an explicit way of
describing the stabilization map.

For each finite pointed set $X_+$ the chain-complex $SQ'_*(X_+)$
in degree $n$   is the free  $k$-module generated by all
families ${\chi}(\varepsilon_1, \ldots, \varepsilon_n)$ of
pairwise disjoint subsets of $X$ indexed by $n$-tupels of elements
$\varepsilon_i \in \{0,1\}$. We divide out all elements that map a
face or a diagonal of  the cube to the empty set. The result of
this normalization process is $SQ_*(X_+)$. The boundary is
$\delta := \sum_{i=1}^n (-1)^i(P_i - R_i - S_i)$. Here
$$ \begin{array}{rcl}
R_i(\chi)(\varepsilon_1, \ldots,
\varepsilon_{n-1}) & = & \chi(\varepsilon_1, \ldots,
\varepsilon_{i-1}, 0, \varepsilon_i, \ldots, \varepsilon_{n-1}), \\
S_i(\chi)(\varepsilon_1, \ldots, \varepsilon_{n-1}) & = & \chi(\varepsilon_1,
\ldots, \varepsilon_{i-1}, 1, \varepsilon_i, \ldots,
\varepsilon_{n-1})
\end{array}$$
and  $P_i(\chi)$ is given by the pointwise union of $R_i(\chi)$ and $S_i(\chi)$.

\Def{
On a generator $f\colon [m] \ra [n+1]$ with $n \geq 1$
the stabilization map $\stab\colon  \ga_{n+1} \ra SQ_n$ is defined
as $\stab(f) := \chi(f)$ where
$\chi(f)(\varepsilon_1, \ldots, \varepsilon_{n})$ is the empty set
for all $n$-tuples which are not of the form
$(0,\ldots,0,\underbrace{1,\ldots,1}_i)$ for  $0 \leq i \leq n$. On
these  tuples the value of $\chi(f)$ is the
preimage of \, $i+1$ under the map $f$:
$$\chi(f)(0,\ldots,0,\underbrace{1,\ldots,1}_i) := f^{-1}(i+1).$$
For $n=0$ we use the convention that $\chi(f)() = f^{-1}(1)$. }

\Exa{Let $f$ be the following map of pointed sets

$$
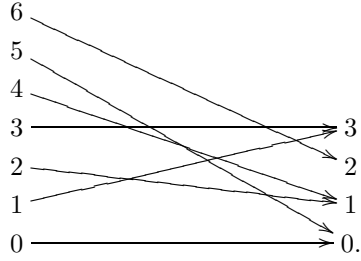
\xymatrix@R=0.5ex{
{6} \ar[rrrrdddd]& & & & {}\\
{5} \ar[rrrrddddd]& & & & {}\\
{4} \ar[rrrrddd]& & & & {}\\
{3} \ar[rrrr]& & & & {3}\\
{2} \ar[rrrrd]& & & & {2}\\
{1} \ar[rrrruu]& & & & {1}\\
{0} \ar[rrrr]& & & & {0.} }$$ Then $\chi(f) \in SQ_2$ is the cube
$$ \left(\begin{matrix}
\{2,4 \} & \{6\} \\
\varnothing & \{1,3 \} \end{matrix}\right). $$}

We have to justify that the map $\stab$ deserves the name
`stabilization'. First we will
show that its image is a subcomplex in the cubical construction.

\Lem{The stabilization induces a map of chain complexes
\label{boundary} from $G(\s^1)_*$ to $Q_{*-1}(G).$ }
\begin{proof}
The boundary of $\stab (f)$ with $f\colon [m] \ra [n+1]$ is given as $\delta
(\chi(f)) = \sum_{i=1}^{n} (-1)^i(P_i - R_i - S_i)(\chi (f))$ As
$\chi (f)$ gives the empty set on every $n$-tuple $(\varepsilon_1,
\ldots, \varepsilon_n)$ which is not of the form $(0, \ldots,0, 1,
\ldots,1)$, the summands $R_i(\chi (f))$ and $S_i(\chi (f))$ are
degenerate except for  $R_0(\chi (f))$ which equals $\chi(d_{n+1}(f))$
and $S_n\chi (f)$ which corresponds to $\chi(d_0(f))$. The
summands $P_i\chi (f)$
 evaluated on an $n$-tuple $(0,\ldots,0,\underbrace{1,\ldots,1}_j)$ give
$\chi (f)(0,\ldots,0,\underbrace{1,\ldots,1}_j)$ for $i \le n-j$,
$\chi(f)(0,\ldots,0,\underbrace{1,\ldots,1}_{j+1})$  for $i>n-j+1$ and
the union of $\chi(f)(0,\ldots,0,\underbrace{1,\ldots,1}_j)$ and $\chi
(f)(0,\ldots,0,\underbrace{1,\ldots,1}_{j+1})$  for $i=n-j+1$.
and these are exactly the values of the face maps $d_{n-i+1}$ for
$i=1,\ldots,n$.
\end{proof}

Thus the stabilization induces a well-defined map $\stab\colon 
\HH_{n+1}(G) \ra \HG_n(G)$ for all $n$ greater or equal to zero and
for all $\ga$-modules $G$. The construction of $\stab$ is
analogous to the one in \cite{EM}, where Eilenberg and MacLane
considered the stabilization map from the homology of
Eilenberg-MacLane spaces to the homology of the corresponding
spectrum. Their result gives one concrete example for the
connection of Hochschild homology and Gamma homology via the
stabilization map.

\Ex{ Let $C$ be an abelian group. The $(n+1)$st Hochschild homology
of the group algebra $k[C]$ with coefficients in the ground ring $k$
is nothing but the group homology of $C$ with coefficients in $k$,
i.e., the $k$-homology of the Eilenberg-MacLane space $K(C,1)$. The
stabilization map has Gamma homology of $k[C]$ as its target and
this is the $k$-homology of the Eilenberg-MacLane spectrum $HC$ (see
\cite{RiRo}). In this case the stabilization map $\stab\colon 
H_{n+1}(K(C,1);k) \ra Hk_nHC $ coincides with the one from
\cite[p.547]{EM}.}

As in
\cite[p.546]{EM} we denote the subcomplex generated by the image of the
stabilization map by $SQ_*^{(0)} \otimes_\ga G = Q_*^{(0)}(G)$. In the
following, we will assume that $G$ is reduced.

\Lem{Assume that $G$ is reduced. The homology of the subcomplex
$Q_*^{(0)}(G)$ is precisely
Hoch\-schild homology of $G$ shifted by one. }
\begin{proof}
For reduced functors $G$, we have that $G[n] \cong \ga_n \otimes_\ga
G$ coincides with $G[n]/G[0]$ which is the cokernel of
$$ \ga_0 \otimes_\ga G \lra \ga_n \otimes_\ga G.$$
On the surjective generators of $\ga_n[m]/\ga_0[m]$ the stabilization map is
injective and the surjective maps in $\ga_n[m]$ correspond to the
normalized chains for Hochschild homology.

Assume that the image of an element $y \in G[n+1]$ in
$Q_n^{(0)}(G)$ is a boundary, $\stab(y) = \delta (\rho)$ for one $\rho$ in
$Q_{n+1}^{(0)}(G)$. The proof of Lemma \ref{boundary} shows that the terms in
the boundary of $\rho$ are in one-to-one correspondence with boundary terms in
the preimage of the stabilization map. Therefore there is an element $y' \in
G[n+2]$ with $\stab(y') = \rho$ and with boundary exactly $y$.
\end{proof}

With the help of the explicit shape of the stabilization map we
can now indicate one reason why the periodicity sequence collapses
after stabilization. Note that $\stab\colon \HH_1 \ra \HG_0$ is an
isomorphisms.

\Thm{\label{thm:collapse} For every $\calf$-module $F$ the
composition of the $B$-map with the stabilization map is trivial
for all $n \geq 1$}

$$ \xymatrix{
 {\HC_n(F)} \ar[rd]_{0} \ar[r]^{B} & {\HH_{n+1}(\mu^*F)} \ar[d]^{\stab} \\
  & {\HG_n(\mu^*F)} }$$

\begin{proof}
In order to prove the claim we will actually show more. We claim that
every element $x \in \calf_n \otimes_\calf F$  is sent to a linear
combination of degenerate elements in $SQ_n \otimes_\ga
\mu^*F$ for $n \geq 1$. Without loss of generality we may assume that $x$ is a
generator, \ie, $x = f \otimes y$ with $f \in \calf_n(\unm)$.  By
definition
$$\stab \circ B(x) = \stab((1-\tau) \circ s \circ N)(f \otimes y).$$
The terms $\stab \circ s \circ \tau^i(f) \otimes y$ are
degenerate, because the composition $\stab \circ s \circ
\tau^i(f)(\varepsilon_1, \ldots, \varepsilon_{n-1},0)$ is empty
for all $(\varepsilon_1, \ldots, \varepsilon_{n-1}) \neq
(0,\ldots,0)$ by definition. As precomposition with $\nu$ causes a
shift by one, the map $s$ causes a trivial preimage of $0$ and
therefore we obtain $\varnothing$ as a value on the $n$-tuple
$(0,\ldots,0)$ as well.

The other terms $\tau \circ s \circ \tau^i(f)$ are degenerate
because they give the empty set on $n$-tuples $(\varepsilon_1,
\ldots, \varepsilon_{n})$ with $\varepsilon_{n-1} \neq
\varepsilon_{n}$: if $\varepsilon_{n-1}$ is $1$ and
$\varepsilon_{n}=0$ then this element gives the empty set by
definition of the stabilization map. In the other case, the map
$\tau \circ s$ has an empty preimage of $1$ and thus $\nu^*$ causes
$\varnothing$ as a value on $(0,\ldots,0,1)$. In particular we obtain
a degenerate element in the case $n=1$.
\end{proof}

\Rek{
As the cubical complex $SQ_*$ is a resolution of the functor $t$ and
as we know that $\nu^*$ is exact, we obtain a resolution $\nu^*SQ_*
\ra \nu^*t$. As we have an isomorphism between $\tor_*^\calf (\nu^*t,
F)$ and cyclic Gamma homology from degrees bigger than one, it
suffices to construct a stabilization map from $\HC_n(F)$ to
$\HG_{n-1}(\mu^*F)$ for $n$ bigger than $2$. As the stabilization map vanishes
on the image of $B$ and as the first column of the bicomplex for cyclic
homology gives rise to $\HH_n(\mu^*F)/B(F(\underline{n-1}))$ we define
the cyclic stabilization map, $\stab_C$, as 
$$\stab_C\colon  \HC_n(F) \lra \HGC_{n-1}(F), \quad \stab_C = \stab \circ \pi,
\text{ for } n > 2 $$ where $\pi$ is the projection from $\HC_n(F)$
to $\HH_n(\mu^*F)/B(F(\underline{n-1}))$.

Note, that this does \emph{not} give rise to a well-defined map from
$\HC_2(F)$ to $\HGC_1(F)$ that is compatible with both periodicity
sequences: we would have to compose $\stab \circ \pi$ with the map
from $\HG_1(\mu^*F)$ to $\HGC_1(F)$, but then $\delta$ composed with
that map has to be trivial, but the $S$-map from $\HC_2(F)$ to
$\HC_0(F)$ is non-trivial in general. }

 \mbox{} \pbn

\end{document}